\documentclass[10pt]{amsart}%
\usepackage{amsmath}%
\usepackage{amssymb}%
\usepackage{amscd}%
\usepackage{color}%
\usepackage[matrix, arrow]{xy}
\usepackage[OT2,T1]{fontenc}
\def\ol#1{\overline{#1}}
\def\wh#1{\widehat{#1}}
\theoremstyle{plain}
    \newtheorem{theorem}{Theorem}[subsection]
    \newtheorem{thm}[theorem]{Theorem}
    
    \newtheorem{proposition}[theorem]{Proposition}
    \newtheorem{prop}[theorem]{Proposition}
    \newtheorem{lemma}[theorem]{Lemma}
    \newtheorem{corollary}[theorem]{Corollary}
     \newtheorem{cor}[theorem]{Corollary}
\newtheorem{assumption}[theorem]{Assumption}
\theoremstyle{definition}
    \newtheorem{definition}[theorem]{Definition}
    \newtheorem{defn}[theorem]{Definition}
    
    \newtheorem{example}[theorem]{Example}
    \newtheorem{remark}[theorem]{Remark}
    \newtheorem{rem}[theorem]{Remark}
    
\def\Alphabet{A,B,C,D,E,F,G,H,I,J,K,L,M,N,O,P,Q,R,S,T,U,V,W,X,Y,Z}
\def\alphabet{a,b,c,d,e,f,g,h,i,j,k,l,m,n,o,p,q,r,s,t,u,v,w,x,y,z}
\def\endpiece{xxx}
\def\makeAlphabet[#1]{\expandafter\makeA#1,xxx,}
\def\makealphabet[#1]{\expandafter\makea#1,xxx,}
\def\makeA#1,{\def\temp{#1}\ifx\temp\endpiece\else%
\mkbb{#1}\mkfrak{#1}\mkbf{#1}\mkcal{#1}\expandafter\makeA\fi}%
\def\makea#1,{\def\temp{#1}\ifx\temp\endpiece\else\mkfrak{#1}\mkbf{#1}\expandafter\makea\fi}%
\def\mkbb#1{\expandafter\def\csname bb#1\endcsname{\mathbb{#1}}}
\def\mkfrak#1{\expandafter\def\csname fr#1\endcsname{\mathfrak{#1}}}
\def\mkbf#1{\expandafter\def\csname b#1\endcsname{\mathbf{#1}}}
\def\mkcal#1{\expandafter\def\csname c#1\endcsname{\mathcal{#1}}}
\def\makeop[#1]{\xmakeop#1,xxx,}
\def\mkop#1{\expandafter\def\csname #1\endcsname{{\mathrm{#1}}}} %
\def\xmakeop#1,{\def\temp{#1}\ifx\temp\endpiece\else\mkop{#1}\expandafter\xmakeop\fi}%
\makeAlphabet[\Alphabet]
\makealphabet[\alphabet]
\makeop[Alt,End,Ext,Hom,Sym,Tor]
\makeop[CH,Pic,div]
\makeop[Ker,Im,Coker,Coim,id,pr,isom]
\makeop[Spec,Spf,Spm]
\makeop[Re,Im]%
\makeop[dR,Nis,crys,rig,syn,tor,Cont,mot,cont]
\makeop[Gal,Res,ab,Frob]
\makeop[pol]
\makeop[Var,an,Isoc,univ,Fil]%
\makeop[Gl,Sl]
\makeop[Lie]
\makeop[tay,tr,gr,rk,Tam,tors]
\makeop[Sat]
\makeop[A]
\def\T{T}

\def\isom{\cong}
\def\verk{\circ}
\def\prolim{\varprojlim}
\def\indlim{\varinjlim}

\def\Hom{\operatorname{Hom}}

\def\Tam{\mathrm{Tam}}

\DeclareSymbolFont{cyrletters}{OT2}{wncyr}{m}{n}
\DeclareMathSymbol{\Sha}{\mathalpha}{cyrletters}{"58}
\newcommand{\Q}{\bbQ}
\newcommand{\Z}{\bbZ}
\newcommand{\G}{G}
\newcommand{\F}{\bbF}
\newcommand{\R}{\bbR}

\renewcommand{\A}{\bbA}
\newcommand{\g}{\mathfrak{g}}
\renewcommand{\k}{\mathfrak{k}}
\newcommand{\Gh}{\mathcal{G}}
\newcommand{\Oh}{\cO}
\newcommand{\Dh}{\cD}
\newcommand{\Lh}{\cL}
\newcommand{\Laz}{\cL{az}}
\newcommand{\res}{\mathrm{res}}
\newcommand{\sing}{\mathrm{sing}}
\newcommand{\coh}{\mathrm{coh}}
\newcommand{\nred}{\mathrm{nred}}
\newcommand{\Fr}{\mathrm{Fr}}
\newcommand{\mucoh}{\mu^{\mathrm{coh}}}
\newcommand{\mubk}{\mu^{\mathrm{BK}}}
\newcommand{\mutam}{\mu^{\mathrm{Tam}}}
\newcommand{\MAT}{\mathrm{MAT}}
\newcommand{\BK}{\mathrm{BK}}
\newcommand{\Cl}{\mathrm{Cl}}

\begin{document}
\title{A cohomological Tamagawa number formula} 
\author{Annette Huber and Guido Kings}
\begin{abstract}
	For smooth linear groups schemes over $\bbZ$ we give a cohomological
	interpretation of the local Tamagawa measures as cohomological periods. This
	is in the spirit of the Tamagawa measures for motives defined by Bloch and Kato.
	We show that in the case of tori the cohomological and the motivic Tamagawa measures
	coincide, which reproves the Bloch-Kato conjecture for motives associated to tori.
\end{abstract}
\date{11. August 2009 }
\maketitle

\section*{Introduction}
The Tamagawa number formula for a reductive algebraic group $G$ relates the volume
of $G(\bbA)/G(\bbQ)$ to arithmetical invariants of $G$. Bloch and Kato proposed
a similar formula for all motives, which unfortunately remains largely conjectural.
In their approach, the local Tamagawa measures are given in terms of periods obtained
by comparing different $p$-adic cohomology theories via the Bloch-Kato exponential map.
There is no direct link to
the classical Tamagawa measure for linear groups.

In this article we try to bridge this gap and give a reformulation of
the classical Tamagawa measure. We use the
Lazard isomorphism to define a cohomological "period integral" of
a top class in Lie algebra cohomology over a fundamental class in group 
cohomology of a $p$-adic Lie group. The main result of Section \ref{sect1}
is a Tamagawa number formula for these cohomological periods of reductive groups, see Theorem
\ref{tamglobal}. The key ingredient is an integral version of Lazard's isomorphism proved in
\cite{HKN}.

In the second part of the article, we show that in the case of tori
our cohomological  Tamagawa measures can be identified with the motivic
Tamagawa measures of Bloch and Kato. This gives a proof of the 
Tamagawa number formula for the motives associated to tori over $\bbQ$.
The result as such is not new. By general yoga for
Artin-motives over $\Q$, it suffices to treat the case of $h^0(F)$ where
$F$ is a number field, i.e., the Dedekind $\zeta$-function at $0$ resp. $1$.
This is nothing but the class number formula.

It is our hope that the cohomological Tamagawa number can be related
to the motivic Tamagawa number of Bloch and Kato \cite{BK} for other reductive groups than tori. 
There is some evidence for such a connection.
In the case of
the algebraic group for a division algebra (see Example \ref{exdivision}) the infinite
cohomological period is very closely related to Borel's regulator in \cite{Borel}.
The work of Borel complemented by a comparison of regulators
due to Beilinson and Rapoport provides this relation for the "infinite period", proving the Beilinson conjecture, i.e., the Bloch-Kato conjecture up to a rational factor.
On the other hand, we have shown in \cite{HK2} that the $p$-adic cohomological periods are related to the Bloch-Kato exponential map and the Soul\'e regulator, which are used in the definition of the local motivic Tamagawa measures.

We would like to thank Burt Totaro for a remark on the existence of smooth models of reductive groups.
%

\section{A cohomological Tamagawa measure}\label{sect1}

\subsection{Notation}\label{notations}
We fix the following setting for the whole article. 

Let $\G$ be a smooth linear group scheme over $\Z$ whose generic fibre is connected reductive. 
If $R$ is a ring, we denote $\G(R)$ the $R$-valued points of $\G$ and by $\G_R$ the base
change to $R$.

For any finite place $p$, we choose $\G_p^c$ a  compact open subgroup of $\G(\Q_p)$. 

We say that the pair  $(\G,\G_p^c)$ with $(p\leq\infty)$ has {\em good reduction} at $p$ if $\G_p^c=\G(\Z_p)$ 
and $\G_{\F_p}$ is quasi-split reductive. In particular, $\G^p_c$ is a maximal compact at good reduction primes $p$.

We assume that $(\G,\G_p^c)$ has good reduction at almost all places.

\begin{rem}
We could choose $\G_p^c=\G(\Z_p)$ for all finite places $p$.
Recall also that by \cite{PR} Proposition 3.16, any compact group  of $\G(\Q_p)$ is contained in a maximal compact subgroup.
At the infinite place, we work implicitly with the choice of subgroup 
$\G_\infty^c=\G(\R)$ throughout, which of course is not compact in general.
\end{rem}

Let $\A$ be the ring of adeles of $\Q$ and $\A_f$ the ring of finite adeles.
Note that
\[ \prod_{p<\infty}\G_p^c\subset \G(\A_f) \]
As usual we embed $\G(\Q)$ diagonally into $\G(\A)$. Let
\[ \Gamma=\G(\Q)\cap \prod_{p<\infty}\G_p^c\subset \G(\A_f)\]

\begin{lemma}
\begin{enumerate}
\item
$\G(\Q)\prod_{p <\infty}\G_p^c$ is of finite index in
$\G(\A_f)$.
\item $\Gamma$ is an arithmetic subgroup of $\G(\R)$.
\end{enumerate}
\end{lemma}
\begin{proof}
Note that
\[ \left[\G(\A_f):\G(\Q)\prod_{p <\infty}\G_p^c\right]=
 \left[\G(\A):\G(\Q)\left(\G(\R)\times\prod_{p <\infty}\G_p^c\right)\right]
 \]
Choose an embedding $\G\to \Gl_N$ of the linear group $\G$. Then
$\G(\Q_p)\cap \Gl_N(\Z_p)$ and $\G_p^c$ are commensurable for all $p$.
For almost all $p$, they are even equal. This implies that
$\prod_{p\leq \infty}\G_p^c$ and $\G(\A)\cap \Gl_N(\R\times \hat{\Z})$ are commensurable.
By \cite{PR} Theorem 5.1, $\G(\Q)$ times the latter is of finite index
in $\G(\A)$. 

The description also shows that $\Gamma$ is arithmetic.
\end{proof}

\begin{defn}\label{cgamma}
Let
\[ c_\Gamma=\left[\G(\A_f):\G(\Q)\prod_{p< \infty}\G_p^c\right]=\left[\G(\A)/\G(\Q):\prod_{p\leq\infty}\G_p^c/\Gamma\right]\]
\end{defn}
\begin{example}
Let $K$ be a number field with ring of integers $\Oh$ and consider the Weil restriction $\G=\Res_{\Oh/\Z}\bbG_m$
with $\G_p^c=\G(\Z_p)=(\Oh\otimes\Z_p)^*$ for all primes $p$. 
Then $c_\Gamma$ is the class number of $K$.
\end{example}

\begin{example}Assume that $G_\Q$ is simple, simply connected with $G(\R)$ non-compact. Then by Strong Approximation (\cite{PR} Theorem 7.12) $c_\Gamma=1$.
\end{example}

We need a technical condition.
\begin{assumption}\label{assumption1}
We assume that $\G(\A)/\G(\Q)$ (or equivalently $\G(\R)/\Gamma$) is compact.
\end{assumption}
Note that the condition only depends on $\G_\Q$ and is independent of the choice of $\Gamma$. In fact, by \cite{PR} Theorem 4.12, it is equivalent to $\G_\Q$ being
anisotropic.
\begin{example}\label{extorus}
Let $\T$ be a smooth group scheme over $\Z$ whose generic fibre is a torus. Then
the assumption is equivalent to the condition that the $\Q$-rank of $T$ is zero
(see \cite{PR} Theorem 4.11).
Hence the torus $\Res_{\Oh/\Z}\bbG_m$ does not satisfy the assumption, but the quotient
$\Res_{\Oh/\Z}\bbG_m/\bbG_m$ does. We are going to consider this example in more detail in the second part of the article.
\end{example}

\begin{example} \label{exdivision} Let $D$  be a non-abelian division algebra  over $\Q$ with center $F$. Let
$\Dh$ be a maximal $\Z$-order of $D$. (They exist, see e.g. \cite{CR} Corollary 26.6.)
Let $\G$ be the group scheme of units of a $\Dh$, ie. given by the functor
\[ \G(A)=(\Dh\otimes A)^*\]
for all rings $A$.
The group $\G$ is easily seen to be smooth. $\G_\Q$ is reductive. However,
it does not satisfy our assumption because it contains a diagonal torus $\bbG_m$. 
Let $H$ be the kernel of the reduced norm map. More precisely, for all commutative $\Oh_F$-algebras $A'$ consider
\[ \nred:(\Dh\otimes_{\Oh_F} A')^*\to (A')^*\]
(see \cite{CR} \S 7D, Corollary 26.2).
Let $H'$ be the algebraic group over $\Oh_F$ defined by its kernel. It is a form
of $\Sl_d/F$, in particular simple,
and we let
\[ H=\Res_{\Oh_F/\Z} H'.\]
Then $H_\Q$ is a form of $\Res_{F/\Q}\Sl_d$, in particular semi-simple 
and simply connected. The group is anisotropic and hence it satisfies Assumption \ref{assumption1}. The  smoothness of $H/\Z$ is not obvious. For a direct proof see \cite{Niko} Theorem 2.

The invariant $c_\Gamma$ for  $H$ is $1$ if Strong Approximation  holds for $H'$, i.e., if and only if
$H'(F_v)$ is non-compact for all infinite places $v$ of $F$.
This is the case if the algebra $A'$ is split at the infinite places but not
true in general.
\end{example}

\subsection{The local measure: real case}

Let $\omega$ be an invariant $d$-form on $\G(\R)$. Depending on
the choice of $\omega$, we are going to introduce a measure
on $\G(\R)/\Gamma$. Recall that this quotient was assumed to
be compact. Note that $\G(\R)/\Gamma$ is not a manifold in
general. It may have orbifold singularities.

Note also that $\G(\R)/\Gamma$ is not connected in general.
In order to simplify notation, we will assume $G(\R)$ connected
for the rest of this section.
It is easy to extend the formulae to the general case.

Let $\g$ be the algebraic Lie algebra of $\G_\R$. It can be 
identified with the Lie algebra of the real Lie group
$\G(\R)$.

\subsubsection*{The cohomological measure}
Let $A^i_{\G(\R)}$ be the space of real $C_\infty$-differential forms
on $\G(\R)$.  By definition, $\g^\vee=(A^1_{\G(\R)})^{\G(\R)}$,
the space
of invariant one-forms. Moreover, choose $\Gamma'\subset \Gamma$ a normal subgroup of finite index such that $\G(\R)/\Gamma'$ is a manifold. 
Together with Stoke's Theorem this defines natural maps
\begin{multline*}
\bigwedge^\bullet \g^\vee\to (A^\bullet_{\G(\R)})^{\G(\R)}\subset (A^\bullet_{\G(\R)})^{\Gamma} \\
\to (A^\bullet_{G(\R)/\Gamma'})^{\Gamma/\Gamma'}\to
C^\bullet_\sing(\G(\R)/{\Gamma'},\R)^{\Gamma/\Gamma'}\to C^\bullet_\sing(\G(\R)/\Gamma,\R)
\end{multline*}
The composition is independent of the choice of $\Gamma'$.

\begin{defn}
Let
\[ \cS: H^i(\g,\R)\to H^i_\sing(\G(\R)/\Gamma,\R)\]
the natural map defined by Stoke's theorem, i.e.
 the natural homomorphism induced by the above composition.
\end{defn}


\begin{rem}
Let $K_\infty$ be a maximal compact subgroup of $\G(\R)$ with Lie algebra $\k$.
Recall that the van Est isomorphism identifies relative Lie algebra cohomology
with continuous group cohomology
\[ H^i(\g,\k,\R)\isom H^i_\cont(\G(\R), \bbR).\]
Restricting to $\Gamma\subset \G(\R)$ and identifying the cohomology of $\Gamma$
with $H^i_\sing( K_\infty\setminus \G(\R)/\Gamma,\R)$, we get a map
\[ H^i(\g,\k,\R)\to H^i_\sing( K_\infty\setminus \G(\R)/\Gamma,\R)\]
This map is a crucial ingredient
in Borel's definition of a regulator for number fields.
\end{rem}

\begin{defn}Let $\eta_{\Gamma}\in H^d(\G(\R)/\Gamma,\Q)$ be the fundamental class, i.e.,
the Poincar\'e dual of the cycle $\G(\R)/\Gamma$ in $H_0(\G(\R)/\Gamma,\Q)$.
\end{defn}
 
The fundamental class is a basis of $H^d(\G(\R)/\Gamma,\Q)$.

\begin{rem}If $G(\R)$ is not connected, this construction has to be
carried out in such a way that Lemma \ref{projinfinite} holds. We leave the details to the reader.
\end{rem}

\begin{defn}
Let $c\in H^d(\g,\R)$. Then 
\[ \cS(c)=\pi(\Gamma,c)\eta_{\Gamma}\]
for some real number, the {\em cohomological period at infinity}.
Let $\mucoh_c$ be the unique translation invariant measure on $G(\R)/\Gamma$ normalized by
\[ \mucoh_c(G(\R)/\Gamma)=\pi(\Gamma,c)\ .\]
\end{defn}

\begin{lemma}\label{projinfinite}
Let $\Gamma'\subset\Gamma$ be a normal subgroup of finite index.
Then
\[ \pi(\Gamma',c)=[\Gamma:\Gamma']\pi(\Gamma,c)\]
\end{lemma}
\begin{proof}
This follows directly from the projection formula for the covering $G(\R)/\Gamma'\to G(\R)/\Gamma$.
\end{proof}

Let $\omega\neq 0$ be an invariant $d$-form on $\G(\R)$. It is
a gauge form and hence defines a (non-normalized) 
translation invariant measure $\mutam_\omega$ on the compact homogeneous space $\G(\R)/\Gamma$. 
This is nothing but
the local Tamagawa measure at infinity.

On the other hand, $\omega$ induces a class $[\omega]\in H^d(\g,\R)$.
\begin{prop}\label{taminfinite} 
Let $\omega\neq 0$ be an invariant $d$-form on $\G(\R)$.
Let $\Gamma\subset\G(\R)$ be a discrete subgroup such that $\G(\R)/\Gamma$
is compact. Then
\[ \mutam_\omega=\mucoh_{[\omega]}\]
The cohomological period agrees with the local Tamagawa number.
\end{prop}
\begin{proof} By Lemma \ref{projinfinite} we may replace $\Gamma$ by a normal
subgroup of finite index. We 
choose $\Gamma$ small
enough such that $\G(\R)/\Gamma$ is a manifold. By definition the isomorphism $\cS$
is induced by the Theorem of Stokes for $\G(\R)/\Gamma$, i.e., by
integrating differential forms over cycles.
\end{proof}

\subsection{The local measure: $p$-adic case}
Let $\G$ be a smooth algebraic group scheme over $\Z_p$ of dimension $d$ with connected generic fibre. Let $\g$ be its $\Z_p$-Lie algebra.
Let $\omega$ be an invariant $d$-form on $\G_{\Q_p}$.
Depending on the choice of $\omega$, we are going to introduce a measure on the compact group $\G(\Z_p)$ and compare it
with the classical local Tamagawa measure.

Recall that $\G(\Z_p)$ is a $p$-adic Lie group. 
Let $\Gh\subset G(\Z_p)$ be a compact open subgroup. 
Recall also (\cite{L} V Theorem 2.3.10) that continuous
and locally analytic group cohomology agree for $\Gh$.
We are going to denote it simply by $H^i(\Gh,\Q_p)$.

\begin{prop} The natural restriction map
\[ \res: H^i(\G(\Z_p),\Q_p)\to H^i(\Gh,\Q_p)\]
is an isomorphism.
\end{prop}
\begin{proof}
In \cite{HKN} Theorem 4.3.1 locally analytic group cohomology of
both groups was shown to agree with Lie algebra cohomology.
In particular, the restriction is an isomorphism.
\end{proof}
However, the integral structures on these cohomology groups defined by cohomology with 
coefficients in $\bbZ_p$ differ.

If $\Gh$ is without $p$-torsion, then
by \cite{L} Theorem V 2.5.8 it is a Poincar\'e group with dualizing module $D=\Z_p$.
This means that for every finitely generated $\Z_p$-module $A$ with continuous operation of
$\Gh$, there is a natural isomorphism
\[ H^i(\Gh,A^*)^*\isom H^{d-i}(\Gh,A).\]
(Here $(\cdot)^*$ is the Pontrjagin dual $\Hom(\cdot,\Q_p/\Z_p)$). 
We are interested in the basic case $A=\Z_p$.

\begin{defn}\label{fundclass} Let $\Gh\subset G(\Z_p)$ be an open subgroup
without torsion. We put 
$\eta_\Gh\in H^d(\Gh,\Q_p)$ the image of
$1\in \Z_p=H^0(\Gh,\Q_p/\Z_p)^*$ in $H^d(\Gh,\Z_p)\subset H^d(\Gh,\Q_p)$. It is called {\em fundamental class} of $\Gh$.
\end{defn}
Note that the fundamental class is a generator of $H^d(\Gh,\Z_p)$ as $\Z_p$-module. This condition
determines it up to a factor in $\Z_p^*$.

\begin{rem} There is an alternative point of view. If $\Gh$ is a Poincar\'e group with
dualizing module $D$, then
there is also a natural $\Q_p$-duality. Under the perfect pairing
\[ H^0(\Gh,\Q_p)\times H^d(\Gh,\Q_p)\to \Q_p\]
$\eta_\Gh$ is indeed the dual of $1\in H^0(\G,\Q_p)$.
Details are going to discussed elsewhere (see~\cite{H}).
\end{rem}

\begin{lemma} \label{comp}
Let $\Gh$ be an open torsionfree subgroup of $G(\Z_p)$ without $p$-torsion and
$\Gh'\subset \Gh$ a subgroup of index $N=p^k$. Then
\[ \eta_\Gh|_{\Gh'}= p^k\eta_{\Gh'}\]
\end{lemma}
\begin{proof}This follows from the projection formula for restriction and corestriction.
\end{proof}

Let $c\in H^d(\g,\Q_p)$. We think of it as a volume form on $G(\Z_p)$.
It is represented by some invariant $d$-form $\omega$ on $G(\Z_p)$.
Recall (\cite{L} Theorem V 2.4.9, see also \cite{HKN} Theorem 4.3.1) that the Lazard morphism is an isomorphism
\[ \Laz: H^d(\Gh,\Q_p)\to H^d(\g,\Q_p)\]
Hence $\Laz^{-1}c$ is a multiple of the fundamental class of $G(\Z_p)$. We define $\pi(\Gh,c)\in\Q_p$ by
\[ \Laz^{-1}c= \pi(\Gh,c) \eta_\Gh.\]
\begin{defn} Let $\Gh\subset G(\Z_p)$ be an open subgroup
and $c\in H^d(\Gh,\Q_p)$ a cohomology class.
Choose $\Gh'\subset \Gh$  a torsion free open subgroup.
Let $\mucoh_c$ be the unique Haar measure on $\Gh$ normalized by
\[ \mucoh_c(\Gh')=|\pi(\Gh',c)|_p\]
\end{defn}
\begin{lemma} The cohomological measure is well-defined, i.e., independent of the choice of $\Gh'$.\end{lemma}
\begin{proof} This follows directly from Lemma \ref{comp}.\end{proof}

\subsection{The local Tamagawa measure}\label{local-Tam-section}
An invariant $d$-form $\omega$ defines a local Tamagawa measure $\mutam_\omega$ on $\G_p^c$, see \cite{Weil} 2.2.1. More precisely, let $U\subset G_p^c$,
let $x_1,\dots,x_d$ be local coordinates on $U$ and $\omega=fdx_1\wedge dx_2\wedge\dots\wedge dx_d$ on $U$. Then
\[ \mutam_\omega(U)=\int_{x(U)} |f\circ x^{-1}|_p\]
with respect to the standard measure on $x(U)\subset \Q_p^d$.

\begin{thm}\label{tamfinite} 
Let $\omega\neq 0$ be an invariant $d$-form on $\G_{\Q_p}$. Let $[\omega]$
be the corresponding class in $H^d(\g,\Q_p)$. Then
\[ \mucoh_{[\omega]}=\mutam_\omega\ .\] 
as measures on $G(\Z_p)$.
\end{thm}
\begin{proof}
It suffices to check
\[ \mucoh_{[\omega]}(\Gh)=\mutam_\omega(\Gh)\ .\] 
for some open compact subgroup $\Gh\subset \G(\Q_p)$. 
As $\omega$ is unique
up to scaling, it also suffices to check the assertion for some $\omega$.

Let
\[ \G(p\Z_p):=\ker(\G(\Z_p)\to \G(\F_p))\]
This is a standard group in the sense of Serre, \cite{Serre}, see also the
detailed discussion in \cite{HKN} Section 2.2.
Let $t_1,\dots,t_d$ be coordinates
of the formal group $\hat{\G}$, the completion of $G$ along $e$. Then
\[ t=(t_1,\dots,t_d):\G(p\Z_p)\to \Z_p^d\]
is well-defined with image $(p\Z_p)^d$. The valuation
\[ \omega(x)=\sup_i |t_i(x)|_p\]
turns $\G(p\Z_p)$ into a $p$-valued group. Let $\Gh$ be its saturated
subgroup, described as follows:
In order to unify notation, let $q=p$ for $p\neq 2$ and $q=4$ for $p=2$.
By \cite{HKN} Lemma 2.2.2 
\[ t(\Gh)=(q\Z_p)^d\ .\]
Recall that $\g$ is the $\Z_p$-Lie algebra of the group scheme $\G$. It is a free
$\Z_p$-module of rank $d$. The exponential map induces a homeomorphism
\[ q\g\to \Gh.\]
By \cite{HKN} Example 2.6.8, the $\Z_p$-Lie algebra $q\g$ is nothing
but the integral Lazard Lie algebra $\Lh(\Gh)$ of $\Gh$. We have
\[ H^d(q\g,\Q_p)=\bigwedge^d(q\g)^\vee\isom\Z_p.\]
Let $c$ be a generator.
The main result of \cite{HKN} is the compatibility of the Lazard isomorphism
with integral structures. By loc. cit. Theorem 3.1.1. together with Example 3.3.1, 
we know that
\[ \Laz: H^d(\Gh,\Z_p)\to H^d(q\g,\Z_p)\]
is an isomorphism. In particular, $\Laz^{-1}(c)$ is a generator of $H^d(\Gh,\Z_p)$.
This implies
\[ \mucoh_c(\Gh)=1.\]

We now turn to the Tamagawa measure.
Let $\omega$ be a generator of the space of
invariant algebraic $d$-forms on $\G$. This space is a $\Z_p$-module of 
rank $1$,
hence $\omega$ is well-defined up to a factor in $\Z_p^*$. Restricting to the
cotangent space at $e$, we get a generator of $\bigwedge^d\g^\vee$. 
Hence $q^{- d}\omega$ generates $\bigwedge^d(q\g)^\vee$.
This implies 
\[ c=uq^{- d}[\omega]\in H^d(q\g,\Q_p)\]
with $u\in\Z_p^*$. Without loss of generality, $u=1$.

We now write $\omega$ in coordinates of the formal group. It has the form
\[ \omega=fdt_1\wedge\dots dt_n\]
with $f\in \Z_p[[t_1,\dots,t_n]]$ such that $f(0)$ a $p$-adic unit. Recall that
$t(\Gh)=(q\Z_p)^d$. Hence
\[ |f|_p=|f(0)|_p=1\]
on $t(\Gh)$ and the Tamagawa measure of $\Gh$ with
respect to $c=q^{-d}\omega$ is 
\[ \int_{(q\Z_p)^d} |q^{-d}f|_p=\int_{(q\Z_p)^d}q^d=1.\]
This finishes the proof.
\end{proof}

\begin{rem}The above proof used the existence of a smooth model $G/\Z_p$.
By \cite{BLR} 7.1 Theorem 5 this is not a restriction.
\end{rem}

\begin{cor}\label{finitepoints} If $\omega$ is a $\Z_p$-generator of 
$H^d(G(\Z_p),\Z_p)$, then
\[ \mucoh_{[\omega]}(\G(\Z_p))=p^{-d}|\G(\F_p)|\]
\end{cor}
\begin{proof}
This is true for the local Tamagawa number $\tau_\omega$, \cite{Weil} Theorem 2.2.5.
It also follows directly from the computation in the proof of Proposition \ref{tamfinite}.
\end{proof}

\subsection{The global formula}
Let $\G/\Z$ and $\Gamma$ be as fixed in Section \ref{notations}. Recall the
index $c_\Gamma$ from Definition \ref{cgamma}.
Let $S$ be a finite set of places
including $\infty$ and all primes of bad reduction of our data.

We want to define a cohomological Tamagawa number as a product of
the local Tamagawa numbers for all places. However, the product does
not converge in general. We have to introduce convergence factors.

Let $M$ be the motive of (the quasi-split inner form of) $G$ in the sense
of \cite{Gross} Chapter 1. It has good reduction at $p\notin S$. 

\begin{rem}The motive of the variety $G_\Q$ (say in Voevodsky's triangulated category of motives) is a direct sum of Artin-Tate motives. It has a structure of
Hopf object. In can be shown that $M$ in the sense of Gross 
is the primitive part of $M(G)$ and independent of the choice of inner form.
Details will be discussed elsewhere.
\end{rem}

\begin{defn} For $p\notin S$, let 
\[ E_p(\G,s)=\det(1-\Fr_p^{-1} p^{-s}|M) \]
 be the Euler factor of $M$ at $p$. The partial $L$-function of $G$ is defined
as the Euler product
\[ L_S(\G,s)=\prod_{p\notin S} E_p(\G,s)^{-1}\]
\end{defn}
This is $L_S(M,s)$ in the notation of \cite{Gross} Section 9.
Under Assumption \ref{assumption1} the function has an analytic continuation (with no
pole) to $s=1$ (see \cite{Gross} Proposition 9.4).

\begin{lemma}\label{euler}
 For almost all $p$
\[ E_p(M,1)=p^{-d}|G(\F_p)|\]
\end{lemma}
\begin{proof}
For almost all $p$, $G_{\F_p}$ is quasi-split reductive. In this case
use \cite{Gross} Section~3.
\end{proof}

\begin{example}
Let $K$ be a number field with ring of integers $\Oh$.
If $\G=\Res_{\Oh/\Z}\bbG_m$ is as in Example \ref{extorus}, then
$L_S(M,s)=\zeta_{K,S}(s)$ is the partial Dedekind $\zeta$-function. For
$\G=\Res_{\Oh/\Z}\bbG_m/\bbG_m$, it is $\zeta_{K,S}(s)/\zeta_S(s)$.
In particular, it is regular at $s=1$.
If $G=H$ is as in Example \ref{exdivision} (kernel of the reduced norm of a central simple algebra over a number field $F$), then it is an inner form of
$\mathrm{SL}_{n,F}$ and hence
\[ L_S(G,s)= L_S(\mathrm{SL}_{n,F},s)=\prod_{i=1}^{n-1} \zeta_{F,S}(s+i).\]
Again it is regular at $s=1$.
\end{example}

\begin{defn} \label{taucoh}
Let $c\in H^d(\Lie(\G),\Q)$.
Then the {\em cohomological Tamagawa number} is defined as
\[ \tau^{\coh}(\G,\Gamma)=L_S(\G,1)^{-1} \cdot 
        \prod_{p\notin S} \left(E_p(\G,1)^{-1}\mucoh_c(\G_p^c)\right) \cdot \prod_{p\in S}\mucoh_c(\G_p^c).
\]
\end{defn}
For almost all $p$,
\[ E_p(\G,1)^{-1}\mucoh_c(\G_p^c)=1\]
by Corollary \ref{finitepoints} and Lemma \ref{euler}. Hence the product
converges. As usual, $\tau^{\coh}(G,\Gamma)$ is independent of the choices of $S$ and $c$.

\begin{rem} If $\G_\Q$ is semi-simple, then the convergence factors are not
necessary. We have
\[ \tau^\coh(\G,\Gamma)=\prod_{p\leq \infty}\mucoh_c(\G_p^c)\]
\end{rem}

Let $\tau^\Tam(\G)$ be the Tamagawa number of the group $\G$ in the sense of Weil, i.e., the volume of $\G(\A)/\G(\Q)$ with respect to the global Tamagawa measure $\mutam$.
\begin{thm} \label{tamglobal}One has the equality
\[ \tau^\Tam(\G)=c_\Gamma \tau^\coh(\G,\Gamma).\]
\end{thm}
\begin{proof}
We have
\[
\tau^\Tam(\G)=\mutam(\G(\A)/\G(\Q))=c_\Gamma \mutam\left( \prod_{p\leq \infty}\G_p^c/\Gamma\right).
\]
We reorganize
\[ \left(\prod_{p\leq \infty}\G_p^c\right)/\Gamma\isom \G(\R)/\Gamma\times\prod_{p<\infty}\G_p^c.\]
By definition of $\mutam$ the assertion now follows from the local identities
Proposition \ref{taminfinite} and Proposition \ref{tamfinite}.
\end{proof}
\begin{rem} \label{smothness}
Rationally, the Lazard isomorphism exists for all choices of $G_p^c$ and hence
the formulation of the theorem is independent of the choice of model $G$. 
By \cite{BLR} 7.1 Theorem 5 all reductive groups $G_\Q$ allow a smooth model
$G$ over $\Z$ in the sense that we need.
\end{rem}

\begin{cor} If $\G_\Q$ is semi-simple and simply-connected, then
\[ 1=c_\Gamma\tau^\coh(\G,\Gamma).\]
If, in addition, $\G^i(\R)$ is non-compact for all simple factors of $\G_\Q$, then
\[ 1=\tau^\coh(\G,\Gamma).\]
\end{cor}
\begin{proof} In this case the Tamagawa number $\tau^\Tam(\G)$ equals $1$. (See \cite{PR} Section 5.3 for a discussion of results on Tamagawa numbers.)
Under the extra assumption we have $c_\Gamma=1$ by Strong Approximation 
\cite{PR} Theorem 7.12 (with $S=\{\infty\}$).
\end{proof}
This includes the case of our Example \ref{exdivision}.

\section{Comparison with the Bloch-Kato conjecture in the case of tori }\label{sect2}

\subsection{Notations}

Let $\T$ be an algebraic torus of dimension $d$ over $\bbQ$. For each $\bbQ$-algebra $A$,
we let $\T(A)$ be the group of $A$-rational points. We denote
by 
\begin{align}
	X^*:=\Hom_{\ol\bbQ}(\T,\bbG_m)&&X_*:=\Hom_{\ol\bbQ}(\bbG_m,\T)
\end{align}
the group of characters (resp. cocharacters) of $T$ defined over $\ol\bbQ$. 
For each field $\bbQ\subset K$ we denote by 
\begin{align}
	X^*_K:=\Hom_K(\T\times_\bbQ K,\bbG_m)&&X_{*K}:=\Hom_{K}(\bbG_m,\T\times_\bbQ K)
\end{align}
the  group of characters (resp. cocharacters) defined over $K$.  
Let 
\begin{align}
	r:=\rk X^*_\bbQ&&r_p:=\rk X^*_{\bbQ_p}&&r_\infty=\rk X^*_{\bbR}
\end{align}
be the $\bbQ$-rank, the $\bbQ_p$-rank and the $\bbR$-rank of $\T$ respectively. 
\begin{assumption}\label{assumption}
To avoid problems with the pole 
of the Riemann zeta function at $1$, we will consider only tori of $\bbQ$-rank $0$.
\end{assumption}

We denote the motive of $\T$ by $V:=h_1(\T)$. Considered as an Artin-Tate motive
this is $V=X_*\otimes\bbQ(1)$.
We can recover the $N$-torsion $\T[N](\ol\bbQ)$
from $X_*(1)$, by observing that
$$
	\T[N](\ol\bbQ)=\Hom_{\ol\bbQ}(X^*\otimes\bbZ/N\bbZ,\bbG_m)\isom X_*\otimes\mu_N.
$$
In particular, we can identify the Tate module of $\T$ with
\begin{equation}
	X_*\otimes\bbZ_p(1)\isom \prolim_n\T[p^n](\ol\bbQ).
\end{equation}
The following set of  points of $\T$ play an important role in the identification of
the motivic points of $V$: By our assumption \ref{assumption} we have
$$
	 \T^1(\bbA):=\bigcap_{\chi\in X^*_\bbQ}\ker(||\chi||_\bbA)=\T(\bbA),
$$
which implies that $\T(\bbA)/\T(\bbQ)$ is compact. For each finite place $v$ of $\bbQ$
we define the maximal compact subtorus $\T^c(\bbQ_v)\subset \T(\bbQ_v)$
by
\begin{equation}\label{Tcdefn}
	\T^c_v:=\bigcap_{\chi\in X^*_{\bbQ_v}}\ker(|\chi|_v),
\end{equation}
where $|\cdot|_v$ is the norm on $\bbQ_v$ normalized by $|p|_p=\frac{1}{p}$. We have
\begin{equation}\label{tcp-seq}
	0\to \T^c_v\to \T(\bbQ_v)\to \bbZ^{r_v}\to 0.
\end{equation}
Define
\begin{equation}
	\Gamma:=\T(\bbQ)\cap (\T(\bbR)\times\prod_p{\T^c_p})\subset \T(\bbA).
\end{equation}
Note that by Dirichlet's unit theorem this is of the form
\begin{equation}
	\Gamma=(\T(\bbQ)\cap \prod_p{\T^c_p})\times E,
\end{equation}
where $E$ is a free group of rank $r_\infty$ (recall that $r=0$ by our Assumption \ref{assumption})
and $(\T(\bbQ)\cap\prod_p{\T^c_p})$ is a finite group.

\begin{definition}
	Let $M$ be an abelian group then we denote the $p$-adic completion by 
	$$
		 M^{\wedge p}:=\prolim_n M/p^nM
	$$
	and the pro-finite completion by
	$$
		 M^\wedge:=\prolim_N M/NM.
	$$
\end{definition}
Consider the Kummer sequence for $\T$
$$
	0\to \T[p^n]\to \T\xrightarrow{[p^n]}\T\to 0.
$$
\begin{lemma}\label{Kummer}
	For each field $k\supset\bbQ$ the Kummer sequence induces isomorphisms 
	$$
		\T(k)^{\wedge p}\isom H^1(k, X_*\otimes\bbZ_p(1))
	$$
	and
	$$
		\T(k)^{\wedge}\isom H^1(k, X_*\otimes\wh\bbZ(1))
	$$
\end{lemma}
\begin{proof}
	The Kummer sequence induces 
	$$
	0\to \T(k)/N\T(k)\to H^1(k,\T[N])\to H^1(k,\T)[N]\to 0
	$$
	and taking the inverse limit we have $\prolim_nH^1(k,\T)[p^n]=0$ (resp. $\prolim_NH^1(k,\T)[N]=0$) 
	as the transition maps are multiplication by $p$ and $H^1(k,\T)$ is finite.
\end{proof}
\subsection{$\bbR$-valued points of tori }

We consider the motive $h_1(\T)$. The Betti-realization is
$h_1(\T)_B=X_*\otimes\bbQ(1)$, which contains $X_*(1)$ as a lattice. 
The de Rham realization is $h_1(\T)_\dR=\Lie \T$,
which has $\Fil^0h_1(\T)_\dR=0$. Bloch and Kato \cite{BK} (5.6) consider 
$D_\infty:=h_1(\T)_{\dR,\bbR}=\Lie \T_\bbR$ and define the $\bbR$-valued points
of the motive $h_1(\T)$ by
$$
	A(\bbR)=\left(D_\infty\otimes_\bbR\bbC/(\Fil^0D_\infty\otimes_\bbR\bbC+X_*(1))\right)^+.
$$
\begin{proposition}
The $\bbR$-valued points of $h_1(\T)$ are given by
$$
	A(\bbR)=\T(\bbR)
$$
and via the identification $D_\infty\isom \Lie \T(\bbR)$ the natural map
$$
	D_\infty\to A(\bbR)
$$
is the exponential map.
\end{proposition}
\begin{proof}
In our case $D_\infty/\Fil^0D_\infty\isom \Lie \T(\bbR)$ so that $D_\infty\otimes_\bbR\bbC\isom \Lie\T(\bbC)$.
We have an exact sequence 
$$
	0\to X_*(1)\to X_*\otimes\bbC\to \T(\bbC)\to 0,
$$
where the last map is the exponential. Hence, $X_*\otimes\bbC\isom \Lie\T(\bbC)$ and we get
$$
	A(\bbR)=\left(\Lie\T(\bbC)/X_*(1)\right)^+=(\T(\bbC))^+=\T(\bbR).
$$
\end{proof}

\subsection{$\bbQ_p$-valued points of tori }

The aim of this section is to identify the $\bbQ_p$-valued motivic points.
The $\bbQ_p$-valued motivic points are by definition (see \cite{BK} (5.6))
$$
	A(\bbQ_p):=H^1_f(\bbQ_p,X_*\otimes\wh\bbZ(1)).
$$
\begin{theorem}\label{p-adic-motivic-points}
	The $\bbQ_p$-valued motivic points of the motive $V=h_1(\T)$ are given
	by
	$$
		A(\bbQ_p)=\T^c_p.
	$$
\end{theorem}
As a first step, we identify the torsion subgroup of $A(\bbQ_p)$.
\begin{lemma}\label{torsionidentification}
The torsion subgroup of $A(\bbQ_p)$ coincides with the one of ${\T^c_p}$:
$$
	(\T^c_p)_\tors=\T(\bbQ_p)_\tors=A(\bbQ_p)_\tors.
$$
\end{lemma}
\begin{proof}
	By definition of $H^1_f(\bbQ_p,X_*\otimes\wh\bbZ(1))$ the torsion
	coincides with the torsion in 
	$$
		H^1(\bbQ_p,X_*\otimes\wh\bbZ(1))=\prod_lH^1(\bbQ_p,X_*\otimes\bbZ_l(1)).
	$$  
	The torsion in $H^1(\bbQ_p,X_*\otimes\bbZ_l(1))$ is $H^0(\bbQ_p,X_*(1)\otimes\bbQ_l/\bbZ_l)$
	because we have an exact sequence
	$$
		0\to H^0(\bbQ_p,X_*(1)\otimes\bbQ_l/\bbZ_l)\to H^1(\bbQ_p,X_*\otimes\bbZ_l(1))\to
		H^1(\bbQ_p,X_*\otimes\bbQ_l(1)),
	$$
	where the first $0$ appears as $H^0(\bbQ_p,X_*\otimes\bbQ_l(1))=0$ for weight reasons.
	On the other hand
	$$
		\T[l^\infty](\ol\bbQ)=\Hom_{\ol\bbQ}(X^*,\mu_{l^\infty})=X_*(1)\otimes\bbQ_l/\bbZ_l.
	$$
	This gives $H^0(\bbQ_p,X_*(1)\otimes\bbQ_l/\bbZ_l)=\T[l^\infty](\bbQ_p)$.
	To conclude, note that the exaxt sequence \eqref{tcp-seq}
	implies that $\T^c_p[l^\infty]=\T[l^\infty](\bbQ_p)$.
\end{proof}
\begin{lemma}
	One has for $l\neq p$
	$$
		H^1_f(\bbQ_p,V_l)=0.
	$$
	In particular, $H^1_f(\bbQ_p,X_*\otimes\bbZ_l(1))$ is
	torsion.
\end{lemma}
\begin{proof}
	Let $I_p$ be the inertia group at $p$, then by definition (\cite{BK} (3.7.1))
	$H^1_f(\bbQ_p,V_l)$ is given
	by the cokernel of the map
	$$
		V_l^{I_p}\xrightarrow{1-\Fr_p^{-1}}V_l^{I_p}.
	$$
	But $V=X_*\otimes\bbQ(1)$ is the Tate twist of an Artin motive, so that
	$\Fr_p^{-1}$ does not have $1$ as an eigenvalue. This implies that $1-\Fr_p^{-1}$ is injective 
	hence surjective.
\end{proof}
We have by Lemma \ref{Kummer}
$$
\T(\bbQ_p)^{\wedge p}\otimes\bbQ_p\isom H^1(\bbQ_p, V_p).
$$
In order to identify $H^1_f(\bbQ_p, V_p)$ in $ \T(\bbQ_p)^{\wedge p}\otimes_{\bbZ_p}\bbQ_p$ we need:
\begin{lemma}\label{T^cidentification}
Denote by $v_p:\bbQ_p^*\to \bbZ$ the $p$-adic valuation, then
one has an exact sequence
\begin{equation}\label{eq9}
	0\to {\T^c_p}^{\wedge p}\to \T(\bbQ_p)^{\wedge p}\to \bbZ_p^{r_p}\to 0.
\end{equation}
In particular,
$$
	{\T^c_p}^{\wedge p}\otimes_{\bbZ_p}\bbQ_p=\bigcap_{\chi\in 
		X^*_{\bbQ_p}}\ker\left(\T(\bbQ_p)^{\wedge p}\otimes_{\bbZ_p}\bbQ_p\xrightarrow{v_p\verk\chi}\bbQ_p \right).
$$
\end{lemma}
\begin{proof}
Consider the exact sequence
$$
	0\to {\T^c_p}\to \T(\bbQ_p)\to \bbZ^{r_p}\to 0.
$$ 
As $\bbZ^{r_p}$ is free, this sequence splits as a sequence of abelian groups, hence
taking the $p$-adic completion is exact. This gives the sequence
\ref{eq9}.
Tensoring this with $\bbQ_p$ over $\bbZ_p$ we get
$$
	0\to ({\T^c_p})^{\wedge p}\otimes_{\bbZ_p}\bbQ_p\to \T(\bbQ_p)^{\wedge p}\otimes_{\bbZ_p}\bbQ_p\to \bbQ_p^{r_p}\to 0.
$$ 
This implies the result.
\end{proof}

The next lemma identifies $H^1_f(\bbQ_p, V_p)\subset \T(\bbQ_p)^{\wedge p}\otimes_{\bbZ_p}\bbQ_p$.
\begin{lemma}\label{h^1_fandlocalpoints}
	Denote by $v_p:\bbQ_p^*\to \bbZ$ the $p$-adic valuation, then
	$$
		H^1_f(\bbQ_p, V_p)=\bigcap_{\chi\in 
		X^*_{\bbQ_p}}\ker\left(\T(\bbQ_p)^{\wedge p}\otimes\bbQ_p\xrightarrow{v_p\verk\chi}\bbQ_p \right).
	$$
In particular, 
$$
	({\T^c_p})^{\wedge p}\otimes_{\bbZ_p}\bbQ_p=H^1_f(\bbQ_p, V_p).
$$
\end{lemma}
\begin{proof}
	Every $\chi\in X^*_{\bbQ_p}$ defines a map
	$$
		H^1(\bbQ_p, V_p)\xrightarrow{\chi}H^1(\bbQ_p, \bbQ_p(1)).
	$$
	By Lemma \ref{Kummer} for $T=\bbG_m$ one has
	$(\bbQ_p^*)^{\wedge p}\otimes_{\bbZ_p}\bbQ_p\isom H^1(\bbQ_p, \bbQ_p(1))$.
	By \cite{HK} Lemma A.1 and Corollary A.2 one has 
	$$
		(\bbZ_p^*)^{\wedge p}\otimes_{\bbZ_p} \bbQ_p\isom H^1_f(\bbQ_p, \bbQ_p(1))
	$$
	and that $H^1_f(\bbQ_p, \bbQ_p(1))$ is the kernel of the valuation map
	$$
		(\bbQ_p^*)^{\wedge p}\otimes\bbQ_p\xrightarrow{v_p}\bbQ_p.
	$$
	On the other hand, $\chi$ defines $H^1_f(\bbQ_p, V_p)\xrightarrow{\chi}H^1_f(\bbQ_p, \bbQ_p(1))$.
	Putting this information together we obtain
	$$
		H^1_f(\bbQ_p, V_p)\subset \bigcap_{\chi\in 
		X^*_{\bbQ_p}}\ker\left(\T(\bbQ_p)^{\wedge p}\otimes\bbQ_p\xrightarrow{v_p\verk\chi}\bbQ_p \right).
	$$
	By Lemma \ref{T^cidentification} we get 
	$$
		H^1_f(\bbQ_p, V_p)\subset{\T^c_p}\otimes\bbQ_p.
	$$
	To show equality, we consider the dimension of both sides. Via the Bloch-Kato
	exponential $D_\dR(V_p)\isom H^1_f(\bbQ_p, V_p)$ we see that the $\bbQ_p$-dimension
	of $H^1_f(\bbQ_p, V_p)$ is $\dim_{\bbQ_p}V_p$. On the other hand, the Euler characteristic
	formula gives
	$$
		\sum_{i=0}^2\dim_{\bbQ_p}H^i(\bbQ_p, V_p)=-\dim_{\bbQ_p}V_p,
	$$
	which implies that 
	$$
		\dim_{\bbQ_p}H^1(\bbQ_p, V_p)=\dim_{\bbQ_p}V_p+\dim_{\bbQ_p}H^0(\bbQ_p,V_p)=\dim_{\bbQ_p}V_p
		+r_p.
	$$
	With the identification $\T(\bbQ_p)^{\wedge p}\otimes_{\bbZ_p}\bbQ_p\isom H^1(\bbQ_p, V_p)$
	and the exact sequence \eqref{tcp-seq}
	we see that $\dim_{\bbQ_p}({\T^c_p}\otimes_{\bbZ_p}\bbQ_p)=\dim_{\bbQ_p}V_p$. This proves the 
	desired result.
\end{proof}
\begin{lemma}\label{p-adicidentification} Under the identification 
$\T(\bbQ_p)^{\wedge p}\isom H^1(\bbQ_p, X_*\otimes\bbZ_p(1))$ by the Kummer sequence,
the subgroup $H^1_f(\bbQ_p, X_*\otimes\bbZ_p(1))\subset H^1(\bbQ_p, X_*\otimes\bbZ_p(1))$
coincides with $({\T^c_p})^{\wedge p}$.
\end{lemma}
\begin{proof}
Consider the diagram
$$
 \begin{xy}\xymatrix{0\ar[r]&({\T^c_p})^{\wedge p}\ar[r]\ar[d]&\T(\bbQ_p)^{\wedge p}\ar[r]\ar[d]&
 \bbZ_p^{r_p}\ar[r]\ar[d]&0\\
 0\ar[r]&({\T^c_p})^{\wedge p}\otimes_{\bbZ_p}\bbQ_p\ar[r]&\T(\bbQ_p)^{\wedge p}\otimes_{\bbZ_p}\bbQ_p\ar[r]&\bbQ_p^{r_p}\ar[r]&0,}
 \end{xy}
 $$
where both rows are exact by Lemma \ref{T^cidentification}. Using the definition
of $H^1_f(\bbQ_p, X_*\otimes\bbZ_p(1))$ as the pull-back of $H^1(\bbQ_p, X_*\otimes\bbZ_p(1))$
to $H^1_f(\bbQ_p,V_p)$, we get that 
$({\T^c_p})^{\wedge p}\isom H^1_f(\bbQ_p, X_*\otimes\bbZ_p(1))$.
\end{proof}
Finally, we can show Theorem \ref{p-adic-motivic-points}:
\begin{proof}[Proof of Theorem \ref{p-adic-motivic-points}] As both $T_p^c$ and $A(\Q_p)$ contain a subgroup of
the form $\bbZ_p^d$ of finite index, it suffices to show that for all $l$
$$
	({\T^c_p})^{\wedge l}\isom A(\bbQ_p)^{\wedge l}.
$$
For $l\neq p$ we have seen that both sides are torsion and the claim follows
from Lemma \ref{torsionidentification}. For $l=p$ the claim follows from
Lemma \ref{p-adicidentification} and the definition of $A(\bbQ_p)$.
\end{proof}

\subsection{Comparison of the motivic with the local Tamagawa measure}\label{comparison-section}

Let $\omega\neq 0$ be a $\T$-invariant algebraic differential form defined over $\bbQ_v$ of top degree on $\T$.
This form defines the local Tamagawa measure $\mu_\omega^\Tam$ on $\T(\bbQ_v)$ (see Section \ref{local-Tam-section}) for all places $v$ 
with the property that
$$
	\mu_\omega^\Tam({\T^c_p})=p^{-\dim \T}\#\cT(\bbF_p)
$$
for almost all $p$. Here $\cT$ is a smooth model of $\T$ over $\bbZ_p$.

We next explain the motivic measures defined by Bloch and Kato on the local points of the motive $h_1(\T)$.
Choose once for all a rational, top degree translation invariant differential form $\omega $ on $\T$.
This gives a linear form $\omega:\bigwedge^d\Lie \T\to \bbQ$ and we denote by
$\omega^\vee\in  \bigwedge^d\Lie \T$ the dual basis.

By definition \cite{BK} 5.9. the motivic measure $\mu^\BK_\omega$ on $A(\bbR)=\T(\bbR)$
equals $\mu_\omega^\Tam$.

Next consider $A(\bbQ_p)={\T^c_p}$. Here the motivic measure $\mu_{\omega}^\BK$
is the Haar measure on ${\T^c_p}$ normalized as follows:
The Bloch-Kato exponential map induces an isomorphism
$$
	\exp_\BK:\Lie\T_{\bbQ_p}\isom H^1_f(\bbQ_p,X_*\otimes \bbQ_p(1)).
$$
By definition of $A(\bbQ_p)={\T^c_p}$, a subgroup of finite index, say $\cT_1\subset A(\bbQ_p)$
is contained in $H^1_f(\bbQ_p,X_*\otimes \bbQ_p)$. Then 
$$
	\exp_\BK^{-1}(\cT_1)=:\cT\subset \Lie\T_{\bbQ_p}
$$
is a $\bbZ_p$-lattice. Choose a basis $t_1,\ldots, t_d$ of $\cT$, then we normalize $\mu_{p,\omega}^\BK$ by 
$$
	\mu_{p,\omega}^\BK(\cT_1)=|\omega (t_1\wedge\ldots\wedge t_d)|_p.
$$
We need the following information about $\exp_\BK$.
\begin{lemma}\label{BK-exp} The
diagram 
$$
\begin{xy}\xymatrix{\Lie\T_{\bbQ_p}\ar[rr]^/-1.1em/{\exp_\BK}\ar[d]_{\exp}& &H^1_f(\bbQ_p,X_*\otimes \bbQ_p(1))
\ar[d]\\
{\T^c_p}^{\wedge p}\otimes{\bbQ_p}\ar[rr]^/-1.1em/{\mathrm{Kummer}}& &H^1(\bbQ_p,X_*\otimes \bbQ_p(1)).
}
\end{xy}
$$
commutes. Here $\exp $ is the exponential map of ${\T^c_p}$.
\end{lemma}
\begin{proof}
After base change to the splitting field one sees with \cite{BK} Example 3.10.1. 
that $\exp_\BK$ is given by the exponential map.
\end{proof}
\begin{proposition}\label{BK-and-Tam-comparison}
For all $v$ the motivic and the local Tamagawa measures coincide
$$
	\mu_{\omega}^\BK=\mu^\Tam_\omega.
$$
\end{proposition}
\begin{proof}
For $v=\infty$ both measures are defined in the same way and there is nothing to show.
For $v=p$ by Theorem \ref{BK-exp}, we have 
$\exp:\cT\isom \cT_1$. On $\cT_1$ we have the translation invariant form $\omega$ and the pull-back
$\exp^*\omega$ is an invariant form on $\cT$. We claim, that this coincides with $\omega$. 
But $\exp$ induces the identity on the tangent spaces, so that
$0^*\exp^*\omega=e^*\omega$, where $0$ and $e$ are the unit sections of $\cT$ and $\cT_1$.
As $e^*\omega$ is the linear form $\omega:\bigwedge^d\Lie \T\to \bbQ$ the claim follows.
Thus the Tamagawa measure of $\cT_1$ is $\omega (t_1\wedge\ldots\wedge t_d)$,
where $t_1,\ldots, t_d$ is the basis of $\cT$ chosen above. 
\end{proof}

\subsection{Global points of tori }

For any number field $K$ let $\MAT(K)$ be the $\Q$-linear abelian category of mixed Artin-Tate motives over $K$.
\begin{rem}Whereas the category of mixed motives in general is conjectural, the 
subcategory of Artin-Tate motives is well-defined. E.g. as (the opposed category of) the heart of the motivic $t$-structure on the full triangulated subcategory of Voevodsky's category $DM_{\mathrm{gm}}(\Spec K,\Q)$ generated by Artin motives over $K$  and all pure Tate motives $\Q(n)$ for $n\in\Z$.
It contains the homological motive of $T$. Its $h_1$ is given by $X_*\otimes\Q(1)$.

The realization functors attach to all objects of $\MAT(K)$ mixed Hodge structures or
$\Gal(\bar{K}/K)$-modules.
\end{rem}
We put
\[ H^1_{\mot}(K,V)=\Ext^1_{\MAT(K)}(\Q(0),h_1(T)).\]

Recall \cite{BK} page 374 bottom, that the global points $A(\bbQ)$ of the motive $V=h_1(\T)$ are defined
as follows: Let
$$
	H^1_{\mot,f}(\bbQ,V):=\ker\left(H^1_{\mot}(\bbQ,V)\to \prod_{p<\infty} H^1(\bbQ_p,V_p)/H^1_f(\bbQ_p,V_p)\right).
$$
Then $A(\bbQ)\subset H^1_f(\bbQ,X_*\otimes\wh\bbZ(1))$ is the preimage of 
$$
	H^1_{\mot,f}(\bbQ,V)\subset H^1_f(\bbQ, X_*\otimes\wh\bbZ(1))\otimes_\bbZ\bbQ
$$
in $H^1_f(\bbQ,X_*\otimes\wh\bbZ(1))$. 
\begin{theorem}\label{global-motivic-points}
	Recall that $\Gamma=\T(\bbQ)\cap(\T(\bbR)\times  \prod_p{\T^c_p})$, then the global
	points of the motive $V=h_1(\T)$ are
	$$
		A(\bbQ)=\Gamma.
	$$
\end{theorem}
For the proof we need several lemmas:
\begin{lemma}\label{motivic-identification}
	One has 
	\begin{align*}
		H^1_{\mot}(\bbQ,V)&\isom \T(\bbQ)\otimes_\bbZ\bbQ
	\end{align*}
	and
	\begin{align*}
		H^1_{\mot,f}(\bbQ,V)&\isom \Gamma\otimes_\bbZ \bbQ.
	\end{align*}
\end{lemma}
\begin{proof}
	Let $t\in T(\Q)$. By the 
	Abel-Jacobi map the homologically trivial cycle $[t]-[1]$ induces an element
	 of $\Ext^1_\MAT(\Q(0),h_1(T))$. This defines a natural map 
	$T(\Q)\otimes_\bbZ\Q\to H^1_{\mot}(\Q,V)$. 
	Let $K/\bbQ$ be the splitting field of $\T$, then $T(K)=(K^*)^{d}$ and hence 
	$$
		H^1_{\mot}(K,V)\isom \T(K)\otimes_\bbZ\bbQ.
	$$
	Using $H^1_{\mot}(\bbQ,V)=H^1_{\mot}(K,V)^{\Gal(K/\bbQ)}$ the result 
	for $H^1_{\mot}(\bbQ,V)$ follows. With this result and Lemma 
	\ref{h^1_fandlocalpoints}
	we have a cartesian diagram
	\begin{equation}\label{mot-cart-diag}
\begin{xy}\xymatrix{
	H^1_{\mot,f}(\bbQ,V)\ar[d]\ar[r]&\prod_p({\T^c_p})^{\wedge p}\otimes_{\bbZ_p}\bbQ_p\ar[d]\\
	\T(\bbQ)\otimes\bbQ\ar[r]&\prod_p\T(\bbQ_p)^{\wedge p}\otimes_{\bbZ_p}\bbQ_p.}
	\end{xy}
	\end{equation}
On the other hand, consider the exact sequence
\begin{equation}\label{gammaTsequence}
	0\to \Gamma\to \T(\bbQ)\xrightarrow{\prod v_p} 
	\prod_p\T(\bbQ_p)^{\wedge p}/({\T^c_p})^{\wedge p}=\prod_p\bbZ^{r_p}_p.
\end{equation}
If we tensor with $\bbQ$ and use $(\prod_p\bbZ^{r_p})\otimes\bbQ\subset \prod_p\bbQ_p^{r_p}$ we get 
$$
0\to \Gamma\otimes_\bbZ \bbQ\to \T(\bbQ)\otimes_\bbZ \bbQ
\xrightarrow{\prod v_p\otimes\id}\prod_p\bbQ^{r_p}_p.
$$
In particular, the diagram \eqref{mot-cart-diag} is cartesian with
$H^1_{\mot,f}(\bbQ,V)$ replaced by $\Gamma\otimes_\bbZ \bbQ$, which proves our claim.
\end{proof}	

\begin{lemma}\label{completedcartesiandiagram}
There is a cartesian diagram 
$$
	\begin{xy}\xymatrix{\Gamma^\wedge\ar[r]\ar[d]&\prod_p{\T^c_p}^{\wedge}\ar[d]\\
 \T(\bbQ)^\wedge\ar[r]&\prod_p\T(\bbQ_p)^{\wedge }.}
 \end{xy}
 $$
\end{lemma}
\begin{proof}
By definition the diagram
$$
	\begin{xy}\xymatrix{\Gamma\ar[r]\ar[d]&\prod_p{\T^c_p}\ar[d]\\
 \T(\bbQ)\ar[r]&\prod_p\T(\bbQ_p).}
 \end{xy}
 $$
is cartesian. Note that the cokernel of
$\prod_p{\T^c_p}\to \prod_p\T(\bbQ_p)$ is $\prod_p\bbZ^{r_p}$ and
hence torsion free. By equation \eqref{gammaTsequence}
the cokernel of $\Gamma\to T(\bbQ)$ is also torsion free. Moreover, as the
$N$-multiplication on a product is the product of the $N$-multiplications, the exact 
sequence
$$
	\prod_p T(\bbQ)\xrightarrow{[N]}\prod_p T(\bbQ)\to \prod_p (T(\bbQ)/NT(\bbQ))\to 0
$$
shows that $(\prod_pT(\bbQ))\otimes\bbZ/N\isom \prod_p(T(\bbQ)\otimes\bbZ/N)$.
This implies that for 
each $N\in \bbZ$ we have a cartesian diagram
 $$
	\begin{xy}\xymatrix{\Gamma\otimes_\bbZ\bbZ/N\ar[r]\ar[d]&\prod_p{\T^c_p}\otimes_\bbZ\bbZ/N\ar[d]\\
 \T(\bbQ)\otimes_\bbZ \bbZ/N\ar[r]&\prod_p\T(\bbQ_p)\otimes_\bbZ\bbZ/N.}
 \end{xy}
 $$
Passing to $\prolim_N$ and observing that this comutes with products and finite fibre 
products, hence with cartesian diagrams, gives the claim.
\end{proof}
\begin{corollary}\label{motivic-profinite}
One has
$$
	H^1_f(\bbQ,X_*\otimes\wh\bbZ(1))=\Gamma^\wedge=\Gamma\otimes_\bbZ\wh\bbZ.
$$
\end{corollary}
\begin{proof}
The first equality follows from Lemma \ref{completedcartesiandiagram}, 
the definition of $H^1_f(\bbQ,X_*\otimes\wh\bbZ(1))$ by the 
cartesian diagram
$$
	\begin{xy}\xymatrix{H^1_f(\bbQ,X_*\otimes\wh\bbZ(1))\ar[r]\ar[d]&
	\prod_p H^1_f(\bbQ_p,X_*\otimes\wh\bbZ(1))\ar[d]\\
	H^1(\bbQ,X_*\otimes\wh\bbZ(1))\ar[r]&\prod_p H^1(\bbQ_p,X_*\otimes\wh\bbZ(1))}
 \end{xy}
 $$
and the identifications in Theorem \ref{p-adic-motivic-points} and Lemma \ref{Kummer}.
The second equality follows from the fact that $\Gamma$ is a finitely generated
abelian group, so that its profinite completion is given by $\Gamma\otimes_\bbZ\wh\bbZ$.
\end{proof}
\begin{proof}[Proof of Theorem \ref{global-motivic-points}] 
The global motivic points are by definition the fibre product of 
$H^1_{\mot,f}(\bbQ,V)=\Gamma\otimes_\bbZ\bbQ$ (Lemma \ref{motivic-identification}) and  $H^1_f(\bbQ,X_*\otimes\wh\bbZ(1))=\Gamma\otimes_\bbZ\wh\bbZ$ (Corollary \ref{motivic-profinite})
over $H^1_f(\bbQ,X_*\otimes\wh\bbZ(1))\otimes\Q=\Gamma\otimes_\bbZ\bbA_f$.
To finish the proof, we have to show that there is a cartesian diagram
$$
	\begin{xy}\xymatrix{\Gamma\ar[r]\ar[d]&
	\Gamma\otimes_\bbZ\wh\bbZ\ar[d]\\
	\Gamma\otimes_\bbZ \bbQ\ar[r]&\Gamma\otimes_\bbZ \bbA_f.}
 \end{xy}
 $$
As $\Gamma$ is a finitely generated abelian group, we can prove this statement 
	for the free part and the torsion part separately. For the free part it follows from
	the standard cartesian diagram
	$$
	\begin{xy}\xymatrix{\bbZ\ar[r]\ar[d]&
	\wh\bbZ\ar[d]\\
	\bbQ\ar[r]&\bbA_f.}
	\end{xy}
	$$
For the torsion part it suffices to note that $\Gamma$ and  $\Gamma\otimes_\bbZ\wh\bbZ$
have obviously the same torsion.
\end{proof}

\subsection{Global invariants}
Recall the definition of the classical Shafarevich group for tori (\cite{Milne} Theorem 9.11)
\[ \Sha(T)=\Ker\left(H^1(\Q,T)\to \prod_{p\leq\infty} H^1(\Q_p,T)\right).\]
On the other hand, Bloch and Kato define (see \cite{BK} Equation (5.13)) 
 a group $\Sha_\BK(M)$ for all motives $M/\Q$.
Using our identifications of local and global points it reads
\[
   \Sha_\BK(h_1(T))=\Ker
      \left(\frac{H^1(\Q,X_*(1)\otimes \Q/\Z)}{\Gamma\otimes \Q/\Z}
     \to \bigoplus_{p\leq \infty} 
       \frac{H^1(\Q_p,X_*(1)\otimes\Q/\Z)}{T_p^c\otimes\Q/\Z}\right),
\]
where by abuse of notation we put $T_\infty^c=T(\R)$.
\begin{defn}The {\em class group} of $T$ is
\[ \Cl(T)=\frac{T(\A_f)}{T(\Q)\prod_{p<\infty}T_p^c}.\]
\end{defn}
\begin{rem}
The order of $\Cl(T)$ is the constant $c_\Gamma$ of Definition \ref{cgamma}.
In particular, the group is finite. If $T=\Res_{K/\Q}\bbG_m$, then
$\Cl(T)$ is the class group of $K$.
\end{rem}

Our aim is to show:
\begin{prop} \label{shavergleich}There is a natural short exact sequence
\[0\to \Cl(T)\to \Sha_\BK(h_1(T))\to \Sha(T)\to 0.\]
\end{prop}
\begin{proof}
Note that
\[ 
	T(\A_f)/\prod_{p<\infty}T_p^c\isom 
		\bigoplus_{p<\infty} T(\Q_p)/T_p^c
\]
(take the direct limit over $S$-adeles for finite sets of places $S$).
By definition of $\Cl(T)$ this implies
\[ 0\to T(\Q)/\Gamma\to \bigoplus_{p<\infty} T(\Q_p)/T_p^c\to \Cl(T)\to 0\]
We abbreviate $I(T)$ for the middle group. Recall that 
$T(\Q_p)/T_p^c\isom \Z^{r_p}$ in our notation. Hence $I(T)$ is torsion free.

Let $n$ be a natural number. We get a commutative diagram with exact rows and
columns
\[ \begin{xy}\xymatrix{
 T(\Q)/\Gamma\ \ar@{^{(}->}[r] \ar@{^{(}->}[d]^{[n]}&I(T)\ar@{->>}[r]\ar@{^{(}->}[d]^{[n]}  & \Cl(T)\ar[d]^{[n]}\\
 T(\Q)/\Gamma\  \ar@{^{(}->}[r] \ar@{->>}[d]&I(T)\ar@{->>}[r] \ar@{->>}[d]  &\Cl(T)\ar[d]\\
  T(\Q)/\Gamma\otimes\Z/n\Z \ar[r]& I(T)\otimes\Z/n\Z\ar@{->>}[r]&\Cl(T)\otimes\Z/n\Z
}\end{xy}\]
By the snake lemma the kernel of the last line is isomorphic to $\Cl(T)[n]$.
This implies that the sequence
\[ 0\to \Cl(T)[n]\to T(\Q)/\Gamma\otimes\Z/n\Z\to \bigoplus_p T(\Q_p)/T_p^c\otimes \Z/n \to \Cl(T)/n\Cl(T)\to 0\]
is exact.
We pass to the direct limit over $n$. Note that
the transition map $\Cl(T)[n]\to \Cl(T)[nm]$ is the natural inclusion. As
$\Cl(T)$ is finite this means
\[ \indlim \Cl(T)[n]=\Cl(T)\]
On the other hand, the transition map $\Cl(T)/n\Cl(T)\to \Cl(T)/nm\Cl(T)$ is
multiplication by $m$. Again by finiteness, this means
\[\indlim \Cl(T)/n\Cl(T)=0\]
We have established the short exact sequence
\begin{equation}\label{sequence}
0\to \Cl(T)\to \frac{T(\Q)\otimes\Q/\Z}{\Gamma\otimes\Q/\Z} \to \bigoplus_p \frac{T(\Q_p)\otimes\Q/\Z}{T_p^c\otimes\Q/\Z}\to 0
\end{equation}
Recall that by abuse of notation $T_\infty^c=T(\R)$, hence the equation
remains valid when the sum runs through $p\leq \infty$.

By Kummer theory for the torus $T$  we have
\[ H^1(k,T(\bar{k}))\isom \frac{H^1(k,T(\bar{k})_\tors)}{T(k)\otimes\Q/\Z}\]
for any field $k\supset\Q$. Hence the defining sequence
for $\Sha(T)$ can be rewritten
\[ 0\to \Sha(T)\to \frac{H^1(\Q,T(\bar{\Q})_\tors)}{T(\Q)\otimes\Q/\Z}\to 
      \prod_{p\leq\infty}\frac{H^1(\Q_p,T(\bar{\Q}_p)_\tors)}{T(\Q_p)\otimes\Q/\Z}.\]
In this sequence we can replace the product by a direct sum because all
global cohomology classes are unramified almost everywhere and unramified
local classes vanish for tori.

Comparing this to the defining sequence of $\Sha_\BK$ yields a commutative
diagram of exact sequences

\[\begin{xy}\xymatrix{
\Sha(T)\ar@{^{(}->}[r] &\frac{H^1(\Q,T(\bar{\Q})_\tors)}{T(\Q)\otimes\Q/\Z}\ar[r]       &\bigoplus\limits_{p\leq \infty}\frac{H^1(\Q_p,T(\bar{\Q}_p)_\tors)}{T(\Q_p)\otimes\Q/\Z}\\
\Sha_\BK(h_1(T))\ar@{^{(}->}[r]\ar[u] &\frac{H^1(\Q,T(\bar{\Q})_\tors)}{\Gamma\otimes\Q/\Z}\ar[r]\ar@{->>}[u]       &\bigoplus\limits_{p\leq\infty}\frac{H^1(\Q_p,T(\bar{\Q}_p)_\tors)}{T_p^c\otimes\Q/\Z}\ar@{->>}[u]\\
\Cl(T)\ar@{^{(}->}[r]\ar[u] &\frac{T(\Q)\otimes \Q/\Z}{\Gamma\otimes\Q/\Z}\ar@{->>}[r]\ar@{^{(}->}[u]       &\bigoplus\limits_{p\leq\infty}\frac{ T(\Q_p)\otimes\Q/\Z}{T_p^c\otimes\Q/\Z},\ar@{^{(}->}[u]
}\end{xy}\]
%
%
where the last line was shown in (\ref{sequence}).
The snake lemma gives the proposition.
\end{proof}
\begin{cor}\label{sha}Let $i(T)$ be Ono's constant \cite{Ono} Section 3.4.
Then
\[ \# \Sha_\BK(h_1(T))= c_\Gamma \cdot i(T)\]
\end{cor}
\begin{proof}
 By definition  $i(T)$ is the order of
\[ \Sha'(T)=\Ker(H^1(K/\Q,T(K))\to H^1(K/\Q,T(\A_K))\]
for big enough $K$ (independent of this choice).
By abuse of notation let $T$ be a model of the torus $T/\Q$ over some
open part of $\Spec \Z$.
By definition
\[ H^1(K/\Q,T(\A_K))=\indlim_S\prod_{v\in S}H^1(K/\Q,T(K_v))\times 
  \prod_{v\notin S}H^1(K/\Q,T(\Oh_v))\]
where the limit is over finite sets of finite places of $K$ and $\Oh_v$ is the ring of
integers of $K_v$.
For any rational prime $p$, we choose a place $v$ of $K$ over $p$. Then
\[ \bigoplus_{v'|p}H^1(K/\Q,T(K_{v'}))\isom H^1(K_v/\Q_p,T(K_v))\]
As in the proof of Proposition \ref{shavergleich} the classes in the image
of $H^1(K/\Q,T(K))$ are unramified almost everywhere and unramifed local classes vanish for tori. Hence we can replace the products by direct sums. Taking the
direct limit over $S$ this means
\[ 0\to \Sha'(T)\to H^1(K/\Q,T(\Q))\to \bigoplus_{p}H^1(K_v/\Q,T(K_v))\]
Passing to the limit over $K$ we get the defining sequence for $\Sha(T)$.
\end{proof}
We also need another global cohomological invariant.
\begin{lemma} \label{globalinv} Let $T_\Q$ be a torus satisfying 
Assumption \ref{assumption1}, i.e, with $\Q$-rank $0$. Then
\[ H^1(\Q,X^*)\isom H^0(\Q,X^*\otimes\Q/\Z).\]
\end{lemma}
\begin{proof}
Under our assumption this follows from the the long exact sequence for
the short exact sequence of discrete Galois modules
\[ 0\to X^*\to X^*\otimes\Q\to X^*\otimes \Q/\Z\to 0.\]
\end{proof}

\subsection{The Tamagawa Number Conjecture for Tori}
Recall that $T/\Q$ is a torus with $\Q$-rank equal to zero.
We now turn to stating the Bloch-Kato conjecture for the motive $h_1(T)=X_*\otimes \Q(1)$. Let $\omega$ be an invariant $d$-form on $T$.
Recall the local measures $\mubk_\omega$ on the local points $A(\Q_p)$ of the motive for $p\leq \infty$ (see Section \ref{comparison-section}). They define a product measure $\mubk$ on
\[ (\prod_{p\leq\infty}A(\Q_p))/A(\Q).\]
This global measure is independent of the choice of $\omega$.

\begin{remark}Note that $h_1(T)$ is pure of weight $-1$. Hence we have to
to use the refined definition \cite{BK} 5.9.1 for the
Tamagawa number. In the same way
as in Definition \ref{taucoh} or in the classical Tamagawa measure 
we introduce convergence factors. 
\end{remark}

The {\em Tamagawa number} of $h_1(T)$ is
\[ \tau^\BK(h_1(T))=\mubk\left((\prod_{p\leq\infty}A(\Q_p))/A(\Q)\right)\]

\begin{thm}\label{TNC}Let $T$ be
a torus of $\Q$-rank equal to $0$. Then the Tamagawa Number Conjecture
of Bloch and Kato (\cite{BK} Conjecture 5.15) holds for the motive $h_1(T)=X_*\otimes\Q(1)$, i.e.,
\[ \tau^\BK(h_1(T))=\frac{\# H^0(\Q,X^*\otimes\Q/\Z))}{\#\Sha_\BK(h_1(T))}.\]
\end{thm}

\begin{proof}[Proof of Theorem \ref{TNC}]
We have already identified local points and local measures, see Proposition \ref{BK-and-Tam-comparison}.
Hence $\mubk=\mutam$.
Using in addition Corollary \ref{sha} and Lemma \ref{globalinv} the claim is equivalent to
\[ \mutam\left(T(\R)\prod_{p<\infty}T_p^c/\Gamma\right)=\frac{\# H^1(\Q,X^*)}{c_\Gamma i(T)}.\]
On the other hand, the classical Tamagawa number formula for tori as proved by Ono (\cite{Ono} Section 5 Main Theorem) reads
\[ \mutam\left(T(\A)/T(\Q)\right)=\frac{\# H^1(\Q,X^*)}{i(T)}.\]
Together with the definition of $c_\Gamma$ this proves the theorem.
\end{proof}
\begin{remark}This clears up a point we had been wondering about: where is
the class number in the Tamagawa number conjecture? 
\end{remark}

\end{document}